\documentclass[
reprint,
superscriptaddress,
amsmath,amssymb,
aps,
pra,
longbibliography,
]{revtex4-2}
\usepackage{float}
\usepackage{xcolor}

\usepackage{graphicx}% Include figure files
\usepackage{dcolumn}% Align table columns on decimal point
\usepackage{bm}% bold math
\usepackage[colorlinks=true, pdfstartview=FitV, linkcolor=blue, citecolor=blue, urlcolor=blue]{hyperref}
\usepackage{dsfont}

\usepackage[separate-uncertainty=true]{siunitx}

\graphicspath{{figures/}}

\begin{document}

\preprint{APS/123-QED}

\title{Activated switching between coexisting limit cycles}

\author{Gabriel Margiani}
\thanks{Authors contributed equally}
\affiliation{Laboratory for Solid State Physics, ETH Z\"{u}rich, CH-8093 Z\"urich, Switzerland.}
\author{Orjan Ameye}
\thanks{Authors contributed equally}
\affiliation{Department of Physics, University of Konstanz, D-78457 Konstanz, Germany.}
\author{Oded Zilberberg}
\affiliation{Department of Physics, University of Konstanz, D-78457 Konstanz, Germany.}
\author{Alexander Eichler}
\affiliation{Laboratory for Solid State Physics, ETH Z\"{u}rich, CH-8093 Z\"urich, Switzerland.}
\affiliation{Quantum Center, ETH Zurich, CH-8093 Zurich, Switzerland}

\begin{abstract}
    Noise-activated switching between coexisting stable states is a fundamental mechanism underlying stochastic dynamics in systems ranging from chemical reactions to neural networks. While this phenomenon is well understood for stationary attractors, it remains largely unexplored for limit cycles, whose periodic motion cannot be described by a static potential landscape. Here we experimentally demonstrate activated switching between two coexisting limit-cycle attractors in a driven nonlinear system of coupled resonators. Specifically, we introduce controlled fluctuations to directly observe the rare stochastic transitions between two limit cycles and measure their dependence on noise intensity and driving strength. The measured switching rates are well described by a large-deviation theory, which replaces the conventional activation barrier by the action along the most probable transition path. Our results extend the concept of activated dynamics from stationary to limit-cycle attractors and establish a framework for modeling stochastic transitions between limit cycles in driven-dissipative systems.
\end{abstract}

\date{\today}

	\maketitle

\begin{figure}[t]
    \includegraphics[width=\columnwidth]{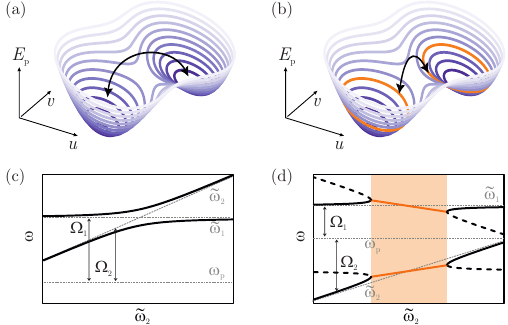}
    \caption{(a)~Schematic illustration of the quasi-potential $E_\mathrm{p}$ of a single KPO in the rotating phase space spanned by $u$ and $v$ (see main text for details). A black arrow indicates stochastic activation between the stable fixed points of $E_\mathrm{p}$, corresponding to the stationary phase states of the KPO. (b)~Visualization of coexisting limit cycles in $E_\mathrm{p}$ (orange lines) and switching between them. (c)~Level repulsion between two dominant advanced signals (solid black lines) around $\Omega_1 = \Omega_2$, as $\tilde{\omega}_2$ is tuned linearly. Idlers on the opposite side of $\omega_\mathrm{p}$ are not presented.
    (d)~Level attraction between an advanced and a delayed dominant signal (solid black lines) and their respective idlers (dashed black lines) around $\Omega_1\simeq-\Omega_2$.
    In the orange shaded range, the modes hybridize (orange lines) and become unstable, leading to limit cycle formation.}
    \label{fig:fig1}
\end{figure}

Most physical models are built around stable equilibrium states, corresponding to local minima of a potential energy landscape. An intuitive picture for an equilibrium state is a ball resting in a valley: the system settles into a stable configuration and stays there until a rare, sufficiently large fluctuation kicks it over a barrier into a neighboring valley. Mathematically, such states are modeled as fixed points, which are the simplest and most common type of attractors of a system's dynamics. For example, a molecule folded into a particular shape can be represented as a fixed point. Any spontaneous change in the molecular shape then corresponds to a noise-activated switch between different fixed points, see Fig.~\ref{fig:fig1}(a). Such noise-activated switching is one of the most fundamental and ubiquitous phenomena in science, underlying stochastic nonlinear dynamics~\cite{aldridge2005noise,Chan_2008,margiani2021fluctuating}, chemical kinetics and gene regulation~\cite{bertram2015mathematical}, neuroscience~\cite{liu2022fixed}, and even game theory and economics~\cite{bertram2015mathematical}.

Not all systems relax to fixed points. Instead, some evolve toward self-sustained, autonomous periodic orbits. These attractors, known as limit cycles, appear in a wide range of contexts, including chemistry~\cite{erban2023chemical} and neuroscience~\cite{jewett1998refinement,leloup1999limit,effenberger2025functional,van1982new}, trained recurrent neural networks~\cite{pals2024trained}, aerodynamic control~\cite{strganac2000identification}, predator–prey models~\cite{may1972limit,pineda2007tale}, and driven–dissipative resonator systems~\cite{aubin2004limit,houri2019limit,lorch2014laser,Czaplewski_2018,dykman2019resonantly,kehrer2025quantum}. In driven–dissipative systems, limit cycles are of particular interest as a resource for generating frequency combs~\cite{kim2020emergence}. Although deterministic in origin, limit cycles are typically difficult to describe analytically, as their properties emerge from strongly nonlinear dynamics rather than from an underlying potential landscape.

Various theoretical models predict the existence of systems with multiple coexisting limit cycles~\cite{galias2022songling,liao2024uniqueness,tang2024bifurcation}. Unlike fixed-point attractors, however, limit cycles are not associated with minima of a potential landscape. Because they are dynamic orbits sustained by a constant balance of driving and dissipation, the familiar picture of noise pushing a system over a simple energy hill breaks down. This raises a fundamental question: what mechanism governs stochastic activation between dynamical attractors, and what determines the transition rate? While indirect evidence of noise-induced transitions between limit cycles has been reported in a granular-gas experiment near a gluing bifurcation~\cite{li2012gluing}, the microscopic switching dynamics could not be controlled or analyzed quantitatively, leaving the mechanism of activation essentially unexplored. 

Here we experimentally realize controllable switching between two coexisting limit-cycle attractors and demonstrate that the measured switching rates match large-deviation theory. To construct these non-equilibrium states, we utilize two coupled electrical resonators, smoothly combining the $\mathbb{Z}_2$ phase-space symmetry of a Kerr parametric oscillator (KPO) with a seeded level-attraction mechanism~\cite{fu2025sideband,Skulte2024realizing} to reliably spawn the coexisting limit cycles. By introducing controlled noise, we induce rare stochastic jumps between these dynamical attractors and systematically extract the transition rates as a function of noise intensity and driving strength. To understand the underlying switching mechanism, we apply a large-deviation approach tailored to limit cycles~\cite{freidlin2012random,e2004minimum,heymann2008geometric,lin2019quasipotential}. By analyzing the stability of points along the orbit, we replace the conventional potential barrier with an action integral calculated along the most probable transition path. Strikingly, the resulting activation cannot be captured by a simple, point-like energy barrier: although the rate retains an Arrhenius-like dependence on noise, its exponent is set by the action along an extended most-probable path through phase space, and it varies with the drive in a way no fixed barrier can reproduce.

To design a system with coexisting limit cycles, we consider two (nearly) identical coupled resonators described by their coordinates $x_1$ and $x_2$
\begin{align}\label{eq:coupled_EOM}
&\ddot{x}_1 + \omega_0^2\left[1 - \lambda \cos(2\omega_\mathrm{p} t)\right]x_1 + \Gamma \dot{x}_1 + \beta x_1^3 - J x_2 = \xi_1\,, \nonumber \\
&\ddot{x}_2 + \omega_0^2 x_2 + \Gamma \dot{x}_2 + \beta x_2^3 -J x_1 = \xi_2\,,
\end{align}
where $\omega_0$ is the resonance frequency, $\Gamma$ is the linear damping rate, $\beta$ is the Duffing nonlinearity, $J$ is the coupling strength, and $\xi_1$ and $\xi_2$ are uncorrelated Gaussian white-noise drives. We parametrically drive only resonator 1 by modulating the quadratic part of its potential at an angular frequency $2\omega_\mathrm{p}\approx 2\omega_0$ with strength $\lambda$. Passing to a frame rotating at $\omega_\mathrm{p}$, we express the coordinates as $x_j(t)=u_j(t)\cos(\omega_\mathrm{p} t)-v_j(t)\sin(\omega_\mathrm{p} t)$ for $j\in \{1,2\}$. We then apply the standard averaging approximation, dropping fast counter-rotating terms to track the slowly varying quadratures $u_j(t)$ and $v_j(t)$~\cite{nonlinearity,Eichler_Zilberberg_book,Seibold_2025,Ameye2025Parametric}.

For $J=0$, resonator 1 acts as a KPO~\cite{Eichler_Zilberberg_book}. Above a critical modulation threshold, the resonator undergoes a $\mathbb{Z}_2$ period-doubling bifurcation and develops two stationary oscillation states with identical amplitudes but with a $\pi$ phase difference. In the rotating-frame phase space spanned by $u$ and $v$, these `phase states' appear as stable fixed points within the system's quasi-potential, see Fig.~\ref{fig:fig1}(a). Small perturbations away from these stationary states evolve with complex eigenfrequencies $\nu_{1,\pm} = \gamma_1 \pm i\Omega_1$. Here, the real part $\gamma_1 < 0$ provides the damping rate that guarantees stability, while the imaginary part $\Omega_1$ denotes the frequency of co- and counter-rotating sidebands in the rotating frame. The latter depend on the drive detuning $\Delta=\omega_\mathrm{p} - \omega_0$ and on amplitude-dependent frequency shifts produced by the Duffing nonlinearity $\beta$~\cite{Huber_2020,heugel2023role,Eichler_Zilberberg_book}. Rotating back to the laboratory frame, interference between the sidebands selects a dominant excitation at $\widetilde{\omega}_1 = \omega_\mathrm{p} + \Omega_1$. Consequently, the sign of $\Omega_1$ dictates in which direction the primary sideband in the lab frame is detuned relative to the pump~\cite{Soriente_2020,Soriente_2021,Dumont_2024}. We call a dominant signal below the pump ($\Omega_j<0$) delayed, whereas a signal above the pump ($\Omega_j>0$) is advanced.

For $J\neq0$, resonator 2 is indirectly driven through resonator 1 and acquires a nonzero stationary response, whose phase is fixed by that of resonator 1. Consequently, the coupled system retains its $\mathbb{Z}_2$ symmetry. In analogy to resonator 1, we denote by $\widetilde{\omega}_2$ the amplitude-renormalized lab-frame eigenfrequency of excitations atop the stationary state of resonator 2. We label the associated rotating-frame eigenvalues as $\nu_{2,\pm}=\gamma_2\pm i\Omega_2$, where $\widetilde{\omega}_2=\omega_\mathrm{p}+\Omega_2$. In a first step, we can now study the possible outcome of what happens when the two resonators share the same sign of $\Omega_j$ and approach one another on the same side of the pump ($\Omega_1\simeq\Omega_2$), see Fig.~\ref{fig:fig1}(c). Their coupling produces level repulsion in analogy to normal-mode splitting between two harmonic oscillators~\cite{Eichler_Zilberberg_book}, resulting in the familiar avoided crossing. 

Conversely, level attraction occurs between a delayed and an advanced signal approaching the pump from opposite sides ($\Omega_1\simeq-\Omega_2$)~\cite{bernier2014unstable,Bernier2018,fu2025sideband}, see Fig.~\ref{fig:fig1}(d). For the parametrically driven responses considered here, we realize this condition using $\beta<0$ and $\Delta<0$: the softening nonlinearity shifts the resonant signal associated with the larger response of resonator 1 below the pump ($\widetilde{\omega}_1<\omega_\mathrm{p}$), where the primary signal becomes delayed. At the same time, the smaller response of resonator 2 remains close to the bare resonance above the pump ($\widetilde{\omega}_2\approx\omega_0>\omega_\mathrm{p}$) and is advanced. As the corresponding sidebands attract, their rotating-frame frequencies approach one another while their decay rates split. Once this splitting overcomes the intrinsic damping, a hybridized eigenvalue pair acquires a positive real part and the stationary phase states lose stability through a Hopf bifurcation~\cite{fu2025sideband,Skulte2024realizing}. In the case realized here, both phase states undergo the Hopf bifurcation at the same point in parameter space, i.e., the phase states give rise to two coexisting limit cycles.

\begin{figure*}[t]
    \includegraphics[width=\textwidth]{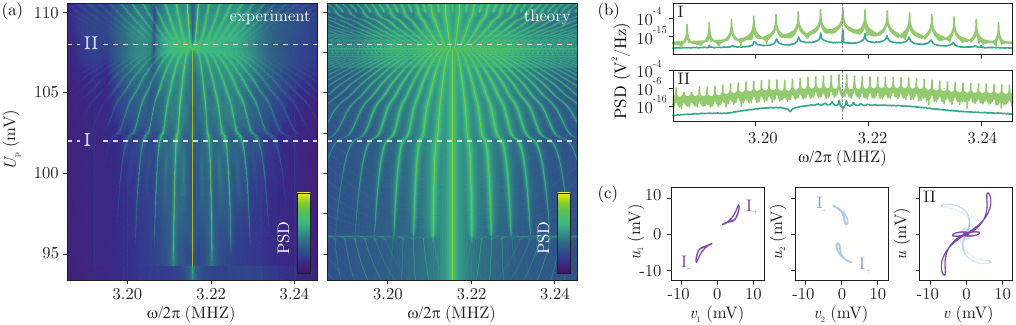}
    \caption{Limit cycle formation. (a)~Measured and simulated oscillation spectrum of resonator 1 around $\omega_\mathrm{p} = \SI{3.2155}{\mega\hertz}$ for varying $U_\mathrm{p}$. (b)~Line cuts of (a) along the lines marked as I and II, with measurement in turquoise and simulation in green. (c)~Measured phase-space trajectories in a frame rotating at $\omega_\mathrm{p}$ for cases I and II. Left panel: the two limit cycles $I_+$ and $I_-$ of resonator 1 in regime I. Middle panel: same as left panel for resonator 2. Right panel: response of resonator 1 and 2 in regime 2. Here and in Fig.~\ref{fig:fig3}, we plot the quadratures of both resonators in a single phase space.}
    \label{fig:fig2}
\end{figure*}

The experimental setup consists of two electrical resistor-inductor-capacitor (RLC) resonators with $\omega_0/(2\pi) = \SI{3222900(100)}{\hertz}$ and $Q=\omega_0/\Gamma = \num{409(1)}$. The resonators are made non-linear and tunable using a variable capacitance diode~\cite{Nosan_2019}, giving a Duffing nonlinearity constant of 
$\beta = \SI{-23(7)e15}{\hertz\squared\per\volt\squared}$. Coupling is achieved inductively via proximity of the two coils~\cite{Heugel_2022}. Placing the coil axes parallel to each other with a separation of roughly \SI{22}{\centi\meter} results in $J = \SI{-1459(4)e9}{\hertz\squared}$. Driving and detection of oscillations at frequency $\omega_\mathrm{p}$ are performed with four separate auxiliary coils that each couple inductively to one of the resonators~\cite{Nosan_2019,margiani2025three}. Cross-inductance between the auxiliary coils is minimized by distance and axis orientation, as well as by the fact that the drive and the relevant oscillation are separated in frequency (at $2\omega_\mathrm{p}$ and $\omega_\mathrm{p}$, respectively). Measurements are performed using a lock-in amplifier (MFLI, Z\"urich Instruments). See SM for details of the setup and calibration procedure.

We drive resonator 1 parametrically as in Eqs.~\eqref{eq:coupled_EOM} with a voltage tone $U_\mathrm{p}\propto\lambda$ at $\omega_\mathrm{p}/(2\pi)= \SI{3.2155}{\mega\hertz}$ and measure the pick-up coil voltages $U_j\propto x_j$, see power spectral density (PSD) of $U_1$ in Fig.~\ref{fig:fig2}(a). For $U_\mathrm{p}<\SI{94}{\milli\volt}$, the drive is too weak to induce parametric oscillation, and we observe only weak fluctuations around $\omega_0$. At $U_\mathrm{p}\approx\SI{94}{\milli\volt}$, resonator 1 undergoes a spontaneous time-translation symmetry breaking in phase space and settles into one of its two phase states~\cite{Eichler_Zilberberg_book}. This symmetry breaking shows up as a dominant peak at $\omega_\mathrm{p}$ in the PSD of $U_1$. The corresponding PSD of $U_2$ is shown in the SM.

For increased $U_\mathrm{p}$, we expect the formation of a limit cycle due to the combination of $\beta<0$ and  $\Delta<0$ as explained above. Indeed, beyond $\SI{94}{\milli\volt}$, we observe the onset of additional, regularly spaced sidebands in the PSD of $U_1$, see regime I in Figs.~\ref{fig:fig2}(a) and (b). To confirm that this is the signature of a limit cycle, we plot the trajectories of $U_1$ and $U_2$ in phase space, see Fig.~\ref{fig:fig2}(c). Indeed, we find that both resonators undergo deterministic and periodic trajectories in $u_i$ and $v_i$. Importantly, the limit cycles emerge from phase state attractors with broken $\mathbb{Z}_2$ symmetry in phase space. As a consequence, there are two different limit cycles separated by a phase difference of $\pi$. The two limit cycles, which we label $I_+$ and $I_-$, form our set of available dynamical attractors, see Fig.~\ref{fig:fig2}(c). We can reach them by switching the parametric drive off and, after a sufficiently long waiting time $\gg 1/\Gamma$, on again; the coupled system then undergoes a new spontaneous symmetry breaking and the limit cycle appears at one of the two positions at random. We thus experimentally establish the coexistence of two different limit cycles in a single system. Note that three co-existing limit cycles were observed in Ref.~\cite{yan2022energy}.

As $U_\mathrm{p}$ approaches $\SI{107.5}{\milli\volt}$, the two limit-cycles expand and merge at the phase-space origin in a so-called homoclinic gluing bifurcation, yielding the characteristic narrow quasi-continuous PSD window~\cite{coullet1984gluing,meron1987gluing,Herrero1998_gluing,li2012gluing}. Beyond this transition (regime II), discrete side peaks reappear but the central peak vanishes, marking a shift between frequency combs with and without a carrier offset~\cite{picque2019frequency}, cf. Figs.~\ref{fig:fig2}(a). In phase space, this absence of a response at $\omega_\mathrm{p}$ reflects a restored symmetry: regardless of initial conditions, the system now traces a single limit cycle centered symmetrically around the origin, see regime II in Fig.~\ref{fig:fig2}(c). The drive amplitude allows us to controllably toggle the system between broken- and restored-symmetry limit cycles~\cite{seitner2017parametric}.

To verify that the observed dynamics correspond to Eqs.~\eqref{eq:coupled_EOM}, we numerically simulate the model using parameters extracted directly from the experiment. The simulations accurately reproduce all primary spectral features, including the limit-cycle onset near $U_\mathrm{p}=\SI{96}{\milli\volt}$, the drive-dependent spacing of the PSD sidebands dictating the cycle periodicity, and the transition from the broken-symmetry regime I to the restored-symmetry regime II at $U_\mathrm{p}=\SI{107.5}{\milli\volt}$ [cf.~Figs.~\ref{fig:fig2}(a) and (b)]. Consequently, we demonstrate a high degree of control over the design and tuning of limit cycles in coupled resonators.

\begin{figure}[t]
\includegraphics[width=\columnwidth]{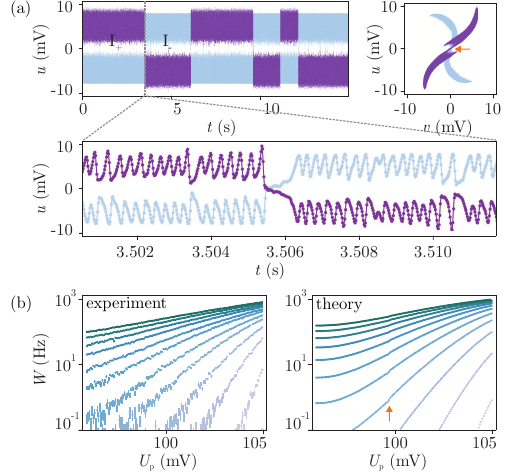}
    \caption{Switching between two limit cycles in regime I with $\sigma = \SI{5}{\milli\volt}$. (a)~Response of resonator 1 (purple) and 2 (blue) as a function of time (left) and in phase space (right), with an orange arrow pointing at the narrow switching channel. Lower line: zoom of one switch. (b)~Switching rate $W$ between limit cycles at $\omega = \SI{3.2155}{\mega\hertz}$ as function of parametric drive $U_\mathrm{p}$ and for noise values between $\sigma = \SI{3.75}{\milli\volt}$ (lightest) and $\sigma = \SI{16.25}{\milli\volt}$ (darkest) in steps of $\Delta\sigma = \SI{1.25}{\milli\volt}$. In the experiment, each data set covers \SI{30}{\second} at a rate of \num{53571} samples per second for each parameter combination. In the simulation, rates were calculated from minimal action $A$ with $C = \SI{1.5(2)}{\kilo\hertz}$ and $b = \SI{7.0(2)e3}{\volt\squared}$ estimated from comparison to the experiment. A vertical arrow marks the discrete jumps in the simulation that we attribute to secondary instabilities.}
    \label{fig:fig3}
\end{figure}

Having established deterministic control in regime I, we drive the system while artificially elevating the noise floor $\xi_j$. We do this by injecting independent white fluctuations of standard deviation $\sigma$ into both resonators using two independent voltage generators. While the system remains confined to a single limit-cycle attractor over short time scales, both $u$ quadratures exhibit irregular sign flips over extended periods, signaling noise-activated switching between $I_+$ and $I_-$, see Fig.~\ref{fig:fig3}(a). We measure the switching rate $W$ as a function of externally applied noise~\cite{margiani2021fluctuating} and find it to scale exponentially with $\sigma^{-2}$ in an Arrhenius-like dependence (see SM). Furthermore, mapping these events in phase space reveals that the transitions occur predominantly along a single narrow channel connecting the two dynamical attractors. This constrained trajectory mirrors the activated switching observed between stationary phase states in a single KPO~\cite{chan2007,chan2008paths,boness2026nonequilibrium}. The confinement of these rare events to a narrow channel isolates a dominant transition path, providing the essential starting point for modeling the transition rates via a minimum-action framework.

After demonstrating stochastic switching between coexisting limit cycles, we formulate these transitions as a noise-activated first-exit process within weak-noise large-deviation theory~\cite{freidlin2012random,e2004minimum,heymann2008geometric,lin2019quasipotential}. In the rotating frame, the averaged equations of motion define a deterministic slow flow where the limit cycles $I_+$ and $I_-$ occupy distinct basins of attraction meeting at a separatrix. The applied fluctuations $\xi_j$ drive rare transitions across this boundary. Consequently, the expected Arrhenius-like switching rate reads $W= C e^{-Ab/\sigma^2}$. Because the limit cycles lack a static potential landscape, the exponent is not a fixed barrier height but the minimum Freidlin-Wentzell action $A$ accumulated along the most probable escape path~\cite{smelyanskiy1997topological,lin2019quasipotential,delacruz2018minimum}.

Beyond the noise dependence, we measure and extract the switching rate as a function of the parametric drive $U_\mathrm{p}$, see Fig.~\ref{fig:fig3}(b). We find that $W$ increases monotonically with the drive. At first sight, this trend appears counterintuitive: increasing $U_\mathrm{p}$ enlarges the amplitudes of the underlying parametric phase states, which naively should suppress switching~\cite{Dykman_1993,Dykman_1998,margiani2021fluctuating,heugel2021ghost,boness2026nonequilibrium}. However, as $U_\mathrm{p}$ increases, the limit cycle deforms and passes progressively closer to the phase-space origin separating the attractors. This geometric deformation shortens the narrow switching channel connecting the limit cycles, strongly reducing the minimum action $A$ required to cross the separatrix. This reduction in the action exponent outweighs the overall increase in oscillation amplitude, yielding the enhanced switching rates.

To validate this mechanism, we numerically compute an approximate value for $A$ for our device parameters using a geometric large-deviation method to find the optimal activation path~\cite{heymann2008geometric, Grafke_2017,heugel2021ghost}. Here, $b$ calibrates the injected noise $\sigma$ against the intrinsic slow-flow fluctuations, and $C$ represents an overall prefactor estimated from the experiment. The resulting theoretical curves reproduce the measured dependence on $U_\mathrm{p}$ [cf.~Fig.~\ref{fig:fig3}(b)]. Since deterministic motion along the cycle costs no action, the corresponding quasipotential remains completely flat over the entire orbit. Escape therefore proceeds via a transverse fluctuation at a preferred phase rather than over a barrier at a single point. Ultimately, limit-cycle activation is governed by the same minimum-action principles as fixed-point escape, extending the Arrhenius paradigm by introducing the geometric requirement of optimizing the departure phase along the dynamical orbit. 

Our numerical analysis reveals additional features as a function of $U_\mathrm{p}$: specifically, the simulated limit-cycle orbits undergo secondary instabilities, including period-doubling bifurcations. These structural transitions manifest as discrete jumps in the theoretical evaluation of the action $A$ in Fig.~\ref{fig:fig3}(b), see orange arrow. While these fine-grained features are not resolved in our current experimental data owing to the limited signal-to-noise ratio, they provide strong motivation for future, more controlled studies. There, investigations should systematically evaluate and compare the activation prefactor $C$ alongside the exponential action~\cite{boness2026nonequilibrium}. Presumably, the signatures of these  dynamical instabilities will appear significantly more pronounced within the prefactor, offering a testbed for the unique signatures of limit cycle activation.

In summary, we have demonstrated controlled, noise-activated switching between two coexisting limit-cycle attractors and shown that the measured rates are captured well by a large-deviation theory. While the switching rates follow an Arrhenius law in noise intensity, they depend on the drive amplitude in a complex way, reflecting the extended shape of the barrier between the limit cycles. This result extends the Arrhenius picture of activated escape, which until now was tied to systems relaxing to potential minima, to attractors that possess no underlying potential at all: in place of a simple energy hill, the switching exponent is set by the action along the most probable escape path, which threads a narrow channel between the two cycles. The rate is thus governed by the global geometry of the periodic orbit rather than by any local property of a single point.

These results establish coexisting limit cycles as controllable, addressable states connected by rare, tunable transitions; the dynamical counterpart of the fixed-point attractors that underlie associative memories and annealing machines in coupled KPOs~\cite{Mahboob_2016, Goto_2016, Dykman_2018, Bello_2019, Okawachi_2020,yamamoto2020coherent,Han_2024,Boehm2025,heugel2021ghost,Margiani_2023,Alvarez_2024,margiani2025three}. This suggests a concrete and, to our knowledge, unexplored direction: building Boltzmann machines and Ising solvers in which the encoding states are limit cycles rather than fixed points, and problem optimization proceeds through activated switching between coexisting orbits. Because the switching rate between limit cycles depends on the drive in a richer, geometry-controlled way than the corresponding rate between fixed points, such machines offer an additional tuning knob over the annealing dynamics, raising the possibility of improved optimization performance. Realizing this will require characterizing not only the action exponent but also the switching prefactor, which for a limit-cycle source is enriched by the neutral phase direction of the orbit and by the cycle period, quantities that do not enter the leading exponential rate~\cite{lin2019quasipotential,smelyanskiy1997topological,boness2026nonequilibrium}. Finally, our approach translates naturally to the quantum regime: the same coupled-KPO ingredients that generate limit cycles here are available in quantum parametric circuits~\cite{Grimm_2019,Puri_2019_PRX,Frattini_2024,beaulieu2025observation,Beaulieu_2025criticality}, where limit-cycle switching may acquire distinct, non-Gaussian signatures. The existence of a controllable classical process strongly motivates the search for its quantum counterpart.

\section*{Acknowledgements}
A.E. and O.Z. acknowledge funding from the Swiss National Science Foundation (SNSF) through the Sinergia Grant No.~CRSII5\_206008/1. O.Z. further acknowledges funding from the Deutsche Forschungsgemeinschaft (DFG) via project numbers 449653034 (Heisenberg), 425217212 (SFB1432), 521530974 (FOR5688), and 545605411 (ANR). We acknowledge the help of Vincent Dumont in proof-reading the manuscript. Claude Opus was used to present the human-produced mathematical results of our work in the supplemental material in a pedagogical and intuitive way.

\bibliography{aipsamp}% Produces the bibliography via BibTeX.

\begin{thebibliography}{75}
\bibitem [{non()}]{nonlinearity}%
  \BibitemOpen
  \href@noop {} {}\bibinfo {note} {In the regime of interest the nonlinearity of mode 2 plays a minor role; when developing intuition one may set the $\alpha x_2^3$ term to zero and treat mode 2 as effectively linear.}\BibitemShut {Stop}%
\end{thebibliography}

\end{document}

\widetext

\begin{center}
  \textbf{\large Supplemental Material:\\
  Activated switching between coexisting limit cycles}

  \vspace{11pt}

  Gabriel Margiani$^{1,\ast}$, Orjan Ameye$^{2,\ast}$, Oded Zilberberg$^2$, and Alexander Eichler$^{1,3}$

  \vspace{11pt}

  \footnotesize
  $^1$\textit{Laboratory for Solid State Physics, ETH Z\"{u}rich, CH-8093 Z\"urich, Switzerland}\\
  $^2$\textit{Department of Physics, University of Konstanz, D-78457 Konstanz, Germany}\\
  $^3$\textit{Quantum Center, ETH Zurich, CH-8093 Zurich, Switzerland}\\
  \vspace{5pt}
  $^\ast$\textit{Authors contributed equally}

\end{center}

\setcounter{equation}{0}
\setcounter{figure}{0}
\setcounter{table}{0}
\setcounter{section}{0}
\setcounter{page}{1}
\setcounter{secnumdepth}{3}
\renewcommand{\theequation}{S\arabic{equation}}
\renewcommand{\thefigure}{S\arabic{figure}}
\renewcommand{\thetable}{S\arabic{table}}
\renewcommand{\thesection}{S\arabic{section}}
\renewcommand{\bibnumfmt}[1]{[S#1]}
\renewcommand{\citenumfont}[1]{S#1}

\maketitle

\section{Device characterization}

\begin{figure}[h]
    \centering
    \includegraphics[width=0.5\linewidth]{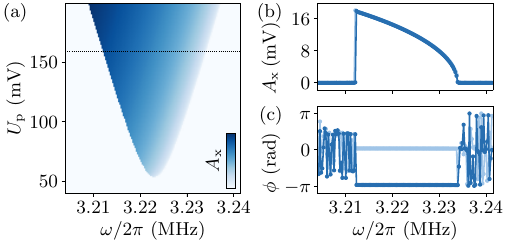}
    \caption{Response of a single resonator to parametric driving. (a) Arnold tongue measured line by line using frequency sweeps from low to high $\omega$. Color linearly encodes the response amplitude from 0 (white) to $A_\mathrm{x} = \SI{20.15}{\milli\volt}$ (dark blue). (b) Amplitude and (c) phase response along the dotted line in (a). The two colors show separate measurements where the resonator responded in each of the two phase states.}
    \label{fig:fig-sm-arnold}
\end{figure}

Our experimental setup consists of two nonlinear RLC resonators that are coupled inductively, cf.~Ref.~\cite{margiani2025three}. To characterize the system, we measure the resonators' linear frequency response as well as the parametric response for a range of $U_\mathrm{d}$. Fig.~\ref{fig:fig-sm-arnold} shows the response of one resonator to a parametric drive while keeping the second resonator far detuned. Repeating the experiment for the other resonator produces very similar data. We extract the resonance frequency $\omega_0/(2\pi)$ and the damping $\Gamma$ from fits to the linear response, while a fit to many parametric frequency sweeps (after applying a correction for three-wave mixing effects~\cite{margiani2025three}) results in values for the parametric threshold $U_\mathrm{th}$ and the Duffing nonlinearity $\beta$. The parametric modulation depth $\lambda$ in Eq.~(1) of the main text is calculated as $\lambda = \frac{2 U_\mathrm{d}}{Q U_\mathrm{th}}$. We, then, tune to the resonators to the same frequency within $\SI{50}{\hertz}\ll \Gamma$ and the same damping within less than one percent~\cite{margiani2025three}. The coupling between the two resonators is calculated from the frequency splitting, which we determine by comparing the bare resonance frequency of the individual resonators to the resonance frequency $\omega_\mathrm{s}$ of the symmetric normal mode appearing when both resonators are tuned to the same resonance frequency. The coupling constant is then approximately given by $J = \omega_0^2 - \omega_\mathrm{s}^2$ for $\omega_0 \gg \omega_0 - \omega_\mathrm{s}$.

\begin{figure}[h]
    \centering
    \includegraphics[width=0.5\linewidth]{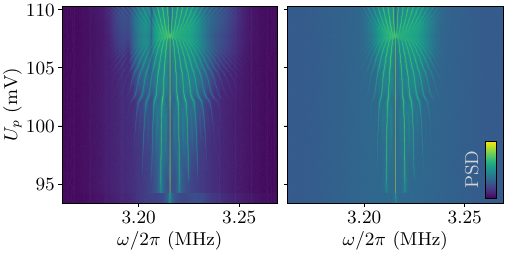}
    \caption{Power spectral density of the driven (left) and undriven (right) resonator for a parametric drive with $\omega_\mathrm{p}/2\pi = \SI{3.2155}{\mega\hertz}$.}
    \label{fig:fig-sm-psd}
\end{figure}

The limit cycle data was measured using a \textit{Zurich Instrumens MFLI}. This lock-in amplifier provides one high resolution input, from which we recorded the response of the driven resonator. To measure the undriven resonator in parallel, we resorted to using an auxilary input with significantly worse noise performance, resulting in noisier data. This can be seen in the power spectral density shown in Fig.~\ref{fig:fig-sm-psd}, with the driven resonator (left panel) showing a lower noise floor (darker background) than the undriven resonators (right panel). The same effect also appears in phase space plots, where the undriven response shows wider lines than the driven response.

Besides the noise, the auxiliary input also features different input delays compared to the main input. We compensate for the resulting phase offset by rotating the response of the undriven resonator by $\delta\phi = \SI{-36}{\degree}$. This value is selected such that the phase response to a linear frequency sweep measured on each of the inputs is the same.

\section{Limit cycle types}
\begin{figure}
    \centering
    \includegraphics[width=\textwidth]{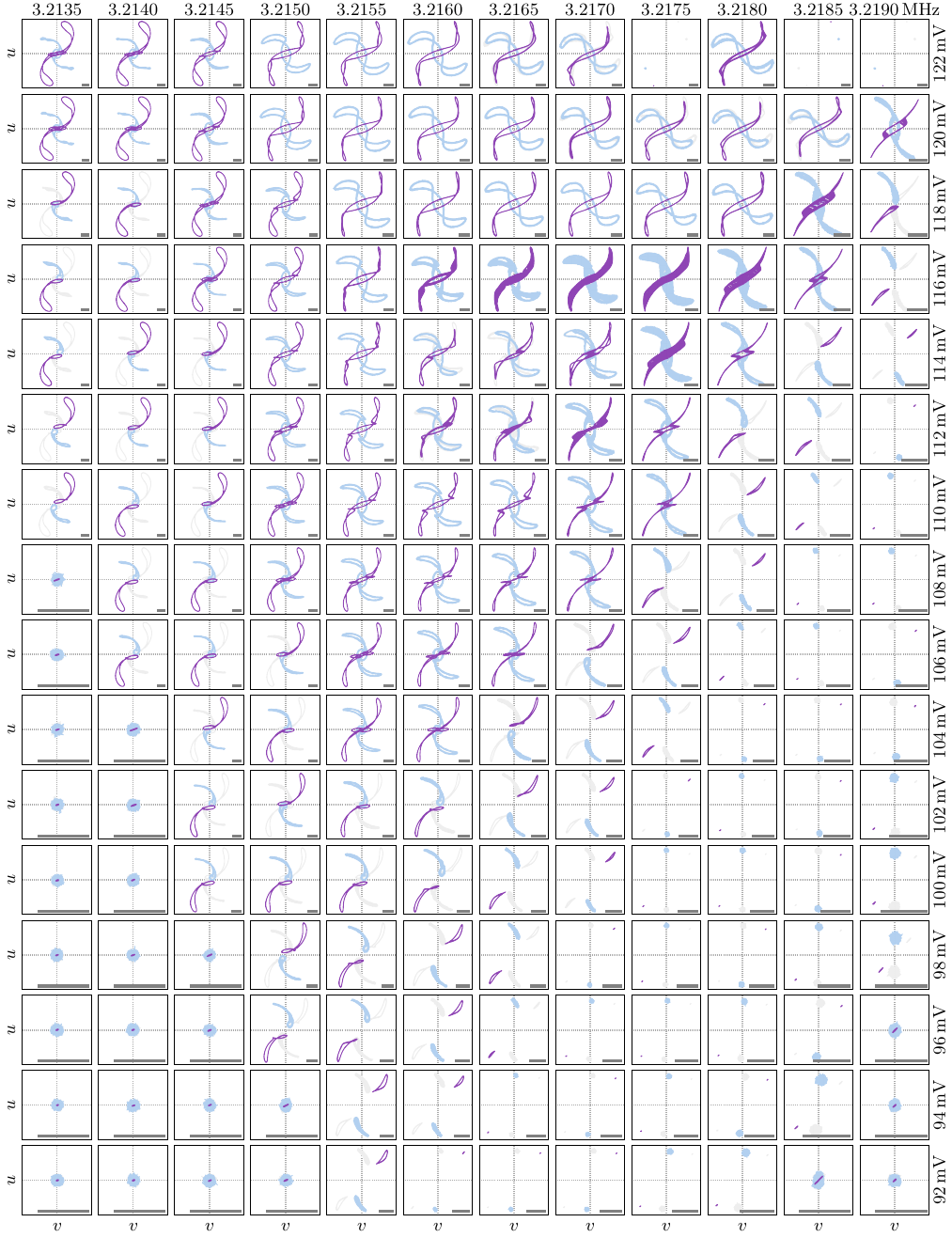}
    \caption{Overview over different limit cycles (see text).}
    \label{fig:fig-sm-table}
\end{figure}

The appearance of the limit cycles depends strongly on the frequency and amplitude of the parametric drive. In Fig.~\ref{fig:fig-sm-table}, we show phase space plots measured at different parametric drives (one per line) and frequencies (one per column), showing different types of limit cycles. Scaling of $u$ and $v$ is given by a gray bar with a length of \SI{3}{\milli\volt} in all panels, while the dotted lines mark the phase space origin. The driven resonator is shown in purple, the undriven one in blue. Light gray lines indicate where a second limit cycle, corresponding to the opposite phase state, would appear (not measured). Most notably we find two limit cycles for example at ($\omega_\mathrm{p}/2\pi = \SI{3.2155}{\mega\hertz}$, $U_\mathrm{p} = \SI{94}{\milli\volt}$) and a single one e.\,g. at ($\omega_\mathrm{p}/2\pi = \SI{3.2155}{\mega\hertz}$, $U_\mathrm{p} = \SI{106}{\milli\volt}$). More chaotic structures can be observed for example at ($\omega_\mathrm{p}/2\pi = \SI{3.2175}{\mega\hertz}$, $U_\mathrm{p} = \SI{116}{\milli\volt}$), while two nearly overlapping cycles show up e.\,g. at ($\omega_\mathrm{p}/2\pi = \SI{3.2175}{\mega\hertz}$, $U_\mathrm{p} = \SI{120}{\milli\volt}$). Furthermore, the differnt noise levels mentioned before are very clearly visible outside of the parametric response, for example at ($\omega_\mathrm{p}/2\pi = \SI{3.2135}{\mega\hertz}$, $U_\mathrm{p} = \SI{92}{\milli\volt}$), while at e.\,g.\ ($\omega_\mathrm{p}/2\pi = \SI{3.2175}{\mega\hertz}$, $U_\mathrm{p} = \SI{92}{\milli\volt}$) the system shows a parametric response without limit cycle.

\section{Switching and noise}

\begin{figure}[h]
    \centering
    \includegraphics[width=0.5\linewidth]{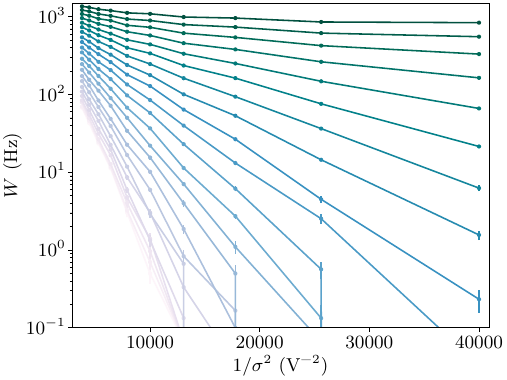}
    \caption{Switching rates as a function of inverse squared noise strength. The plot shows a subset of data points from Fig.~3(b) of the main text. Color encodes $U_\mathrm{p}$ from \SI{93}{\milli\volt} (lightest) to \SI{107.7}{\milli\volt} (darkest).}
    \label{fig:fig-sm-noise}
\end{figure}

To confirm that the noise-dependence of the switching rates follows $W = C_1\exp{-C_2/\sigma^2}$, in Fig~\ref{fig:fig-sm-noise}, we plot the data points from the measurement in Fig.~3(b) of the main text as a function of $1/\sigma^2$. We find good agreement with the exponential dependence for most driving strengths. 
%This might indicate a noise-dependence of the transition paths or limit cycle structure.

\section{Model}

We consider two coupled Duffing resonators:
\begin{align}
\label{eq: supmat eom bare 1}
  \ddot x_1 + \Gamma \dot x_1 + \omega_0^2\left[1-\lambda \cos(2\omega t)\right] x_1 + \beta x_1^3 - J x_2  &= 0 \,,\\
  \ddot x_2 + \Gamma \dot x_2 + \omega_0^2 x_2 + \beta x_2^3 - J x_1 &= 0 \,,
  \label{eq: supmat eom bare 2}
\end{align}
where $x_1$ is parametrically driven at frequency $\omega \approx \omega_0$, while $x_2$ is coupled to $x_1$ but not directly driven.

\subsection{Bare mode}

In the bare mode basis, we define rotating-frame quadratures $u_j$ and $v_j$ via
\begin{equation}
  x_j(t) = u_j(t)\cos(\omega t) - v_j(t)\sin(\omega t), \qquad j = 1, 2,
  \label{eq:bare_quadratures}
\end{equation}
matching the convention used in the main text. After substitution and averaging over the carrier period, $\omega$, the rotating-wave approximation (RWA) yields the slow-flow equations~\cite{Eichler_Zilberberg_book,kovsata2022harmonicbalance}
  \begin{align}
    \dv{u_1}{t} &= \frac{J}{2\omega} v_2
    - \frac{\Gamma}{2} u_1
    + \left(\frac{\omega}{2} - \frac{\omega_0^2}{2\omega} - \frac{\lambda\omega_0^2}{4\omega}\right) v_1
    - \frac{3\beta}{8\omega}(u_1^2 + v_1^2) v_1 \,,\\[1ex]
    \dv{v_1}{t} &= -\frac{J}{2\omega} u_2
    - \frac{\Gamma}{2} v_1
    + \left(\frac{\omega_0^2}{2\omega} - \frac{\omega}{2} - \frac{\lambda\omega_0^2}{4\omega}\right) u_1
    + \frac{3\beta}{8\omega}(u_1^2 + v_1^2) u_1 \,,\\[1ex]
    \dv{u_2}{t} &= \frac{J}{2\omega} v_1
    - \frac{\Gamma}{2} u_2
    + \left(\frac{\omega}{2} - \frac{\omega_0^2}{2\omega}\right) v_2
    - \frac{3\beta}{8\omega}(u_2^2 + v_2^2) v_2 \,,\\[1ex]
    \dv{v_2}{t} &= -\frac{J}{2\omega} u_1
    - \frac{\Gamma}{2} v_2
    + \left(\frac{\omega_0^2}{2\omega} - \frac{\omega}{2}\right) u_2
    + \frac{3\beta}{8\omega}(u_2^2 + v_2^2) u_2 \,.
  \end{align}

It is convenient to package each pair of quadratures into a complex slow envelope $a_j = u_j + i v_j$, in which the four real slow-flow equations above combine into a pair of coupled complex equations,
\begin{align}
  \dot a_1 &= (-\tilde\Gamma + i\delta_1)\,a_1 + i K |a_1|^2 a_1 - i\mu\,a_1^\ast - i g\,a_2\,,\label{eq:model_a1}\\
  \dot a_2 &= (-\tilde\Gamma + i\delta_2)\,a_2 + i K |a_2|^2 a_2 - i g\,a_1\,,\label{eq:model_a2}
\end{align}
in which the physical parameters have been collected as
\begin{equation}
  \tilde\Gamma = \frac{\Gamma}{2},\qquad
  K = \frac{3\beta}{8\omega},\qquad
  g = \frac{J}{2\omega},\qquad
  \mu = \frac{\lambda\omega_0^2}{4\omega},\qquad
  \delta_j = \frac{\omega_j^2 - \omega^2}{2\omega}.
  \label{eq:compact_params}
\end{equation}
In the bare basis, the two resonators share the same natural frequency $\omega_0$, so $\delta_1 = \delta_2 = (\omega_0^2 - \omega^2)/(2\omega)$. We nevertheless retain distinct symbols, for the discussion below. The parametric drive enters only the first resonator, through the term $-i\mu a_1^\ast$, while the linear coupling $g$ acts symmetrically between the two. The compact form~\eqref{eq:model_a1}--\eqref{eq:model_a2} is useful for checking conventions. The lobe physics, however, is most transparent in the normal-mode basis.

\subsection{Normal mode}

When the linear coupling $J$ is comparable to or larger than the detuning and damping scales, the bare-mode basis becomes inconvenient because the linear eigenvectors hybridize the two resonators strongly. It is then more natural to work in the symmetric and antisymmetric normal-mode basis
\begin{equation}
  x_s = \tfrac{1}{2}(x_1 + x_2), \qquad x_a = \tfrac{1}{2}(x_1 - x_2),
\end{equation}
whose linear frequencies $\omega_{s/a} = \sqrt{\omega_0^2 \mp J}$ are the split normal frequencies of the coupled pair. Since $J<0$ in the experiment, $\omega_a<\omega_0<\omega_s$, so the lower-frequency lobe is the antisymmetric normal-mode lobe. Rewriting the equations of motion~\eqref{eq: supmat eom bare 1},~\eqref{eq: supmat eom bare 2} in these coordinates yields
\begin{widetext}
  \begin{align}
    \ddot x_s + \Gamma \dot x_s + \left[\omega_s^2-\frac{\lambda\omega_0^2}{2} \cos(2\omega t)\right] x_s + \beta \left(x_s^2 + 3  x_a^2\right)x_s -\frac{\lambda\omega_0^2}{2} \cos(2\omega t) x_a &= 0 \,,\\
    \ddot x_a + \Gamma \dot x_a + \left[\omega_a^2-\frac{\lambda\omega_0^2}{2} \cos(2\omega t)\right] x_a + \beta \left(x_a^2 + 3 x_s^2\right)x_a -\frac{\lambda\omega_0^2}{2} \cos(2\omega t) x_s &= 0 \,.
  \end{align}
\end{widetext}
In this basis, the parametric drive enters in two distinct ways: It modulates the natural frequency of each normal mode at $2\omega$ (the conventional parametric pump), and, because the original drive sits on resonator $1$ alone rather than on a symmetric superposition, it also generates a time-dependent linear coupling between the symmetric and antisymmetric modes.

Applying the rotating-wave approximation and defining quadratures $u_k, v_k$ for $k=s,a$ via $x_k(t) = u_k(t)\cos(\omega t) - v_k(t)\sin(\omega t)$, we obtain the slow-flow equations
  \begin{align}
  \label{eq: supmat slow normal 1}
    \dv{u_s}{t} &= -\frac{\Gamma}{2} u_s
    - \frac{\omega_s^2 - \omega^2}{2\omega}\, v_s
    - \frac{\lambda\omega_0^2}{8\omega}\,(v_s + v_a)
    - \frac{3\beta}{8\omega}\bigl[(u_s^2 + v_s^2) + (u_a^2 + 3 v_a^2)\bigr] v_s
    - \frac{3\beta}{4\omega}\, u_s u_a v_a \,,\\[1ex]
    \label{eq: supmat slow normal 2}
    \dv{v_s}{t} &= -\frac{\Gamma}{2} v_s
    + \frac{\omega_s^2 - \omega^2}{2\omega}\, u_s
    - \frac{\lambda\omega_0^2}{8\omega}\,(u_s + u_a)
    + \frac{3\beta}{8\omega}\bigl[(u_s^2 + v_s^2) + (3 u_a^2 + v_a^2)\bigr] u_s
    + \frac{3\beta}{4\omega}\, u_a v_s v_a \,,\\[1ex]
    \label{eq: supmat slow normal 3}
    \dv{u_a}{t} &= -\frac{\Gamma}{2} u_a
    - \frac{\omega_a^2 - \omega^2}{2\omega}\, v_a
    - \frac{\lambda\omega_0^2}{8\omega}\,(v_s + v_a)
    - \frac{3\beta}{8\omega}\bigl[(u_a^2 + v_a^2) + (u_s^2 + 3 v_s^2)\bigr] v_a
    - \frac{3\beta}{4\omega}\, u_s u_a v_s \,,\\[1ex]
    \label{eq: supmat slow normal 4}
    \dv{v_a}{t} &= -\frac{\Gamma}{2} v_a
    + \frac{\omega_a^2 - \omega^2}{2\omega}\, u_a
    - \frac{\lambda\omega_0^2}{8\omega}\,(u_s + u_a)
    + \frac{3\beta}{8\omega}\bigl[(u_a^2 + v_a^2) + (3 u_s^2 + v_s^2)\bigr] u_a
    + \frac{3\beta}{4\omega}\, u_s v_s v_a \,.
    \end{align}
After averaging, the time-dependent parametric coupling between the normal modes survives as a static linear coupling of strength $\lambda\omega_0^2/(8\omega)$, mediating both diagonal squeezing and an off-diagonal $s$--$a$ mixing.

In normalized complex envelopes
\begin{equation}
  a_k = \frac{u_k + i v_k}{\sqrt{2}},\qquad
  a_k^* = \frac{u_k - i v_k}{\sqrt{2}},\qquad k = s, a,
\end{equation}
the slow-flow equations take the dissipative-Hamiltonian form
\begin{align}
  \dot a_s = i \frac{\partial H}{\partial a_s^*} - \frac{\Gamma}{2}\, a_s\,, \quad
  \rm{and} \quad
  \dot a_a = i \frac{\partial H}{\partial a_a^*} - \frac{\Gamma}{2}\, a_a\,,
  \label{eq: supmat normal mode dissipative Hamiltonian form}
\end{align}
with the effective Hamiltonian
\begin{align}
  H &= \sum_{k=s,a} \left[\Delta_k\, a_k^* a_k^{\phantom{*}} + K \left(a_k^* a_k^{\phantom{*}}\right)^2
    - \Lambda \left(a_k^* a_k^* + a_k^{\phantom{*}} a_k^{\phantom{*}}\right)
  \right]
  - 2\Lambda \left(a_s^* a_a^* + \text{c.c.}\right)
  + K \left(2\, a_s^* a_s^{\phantom{*}} a_a^* a_a^{\phantom{*}}  + a_s^* a_s^* a_a^{\phantom{*}} a_a^{\phantom{*}} + \text{c.c.} \right)\,,
\end{align}
where $\Delta_k = (\omega_k^2 - \omega^2)/(2\omega)$, $K = 3\beta/(8\omega)$, and $\Lambda = \lambda\omega_0^2/(16\omega)$. Note that the normal-mode envelopes $a_k$ ($k = s, a$) carry an additional $1/\sqrt{2}$ relative to the bare-mode envelopes $a_j = u_j + i v_j$ of Eq.~\eqref{eq:model_a1}; this normalization is what makes $a_k^* a_k^{\phantom{*}}$ the canonically conjugate action coordinate in $H$, so that the equations of motion take the standard Hamiltonian form above. The observed sign $+i\,\partial_{a_k^*} H$ rather than the more familiar $-i\,\partial_{a_k^*} H$ is a consequence of the rotating-frame convention $x_k = u_k\cos\omega t - v_k\sin\omega t$ paired with $a_k = (u_k + iv_k)/\sqrt{2}$; under the alternative phase convention $a_k = (u_k - iv_k)/\sqrt{2}$, the standard sign is recovered.

The Hamiltonian structure also yields a compact form of the normal-mode slow flow equations~\eqref{eq: supmat normal mode dissipative Hamiltonian form},
\begin{align}
  \dot a_a &=(-\tilde\Gamma+i\Delta_a)a_a+i2K(|a_a|^2+|a_s|^2)a_a-i\tfrac{\mu}{2}(a_a^*+a_s^*)+i2K a_s^2a_a^*\,,\label{eq:nm_aa}\\
  \dot a_s &=(-\tilde\Gamma+i\Delta_s)a_s+i2K(|a_a|^2+|a_s|^2)a_s-i\tfrac{\mu}{2}(a_s^*+a_a^*)+i2K a_a^2a_s^*\,,\label{eq:nm_as}
\end{align}
where $\mu/2 = 2\Lambda$. The factor of two relative to the bare-mode pump in Eq.~\eqref{eq:compact_params} reflects the $\sqrt{2}$ rescaling between bare-mode and normal-mode envelopes. These equations are equivalent to the bare quadrature slow flow above, cf.~Eqs.~\eqref{eq:model_a1} and~\eqref{eq:model_a2}. Yet, they make the analytic structure transparent: the linear normal modes diagonalize the conservative mechanics, but the pump is rank-one in normal-mode space. It squeezes $a_a$, squeezes $a_s$, and anomalously couples $a_a$ directly to $a_s^*$, leading to level attraction between sidebands.

\section{Stationary phase states and parametric threshold}
\label{sec:stationary}

We first construct the deterministic skeleton on which the switching problem lives. The skeleton consists of three pieces. First, the origin loses stability and a pair of inversion-related stationary phase states appears. Second, each phase state loses stability through a Hopf bifurcation, producing the two stable limit cycles $I_+$ and $I_-$, cf.~Sec.~\ref{sec:hopf}. Third, the two cycles sit in different basins of attraction. The first two pieces are local deterministic questions and can be treated analytically. The third is global and will enter the large-deviation theory.

We apply the following construction: We use the normal modes to identify the experimentally relevant lobe, but use the passive bare resonator as the ``solvability gauge''. The nonzero stationary states lie on a one-parameter branch labeled by the nonlinear detuning $\vartheta_2$ of the passive resonator. Once this scalar parameter is known, the amplitudes, phases, normal-mode composition, and entries of the linear-response matrix all follow algebraically. The Hopf bifurcation is then located by a scalar condition on the same branch. The remainder of this section builds the stationary branch and locates the parametric threshold; the next section (Sec.~\ref{sec:hopf}) carries the construction through to the Hopf bifurcation and the birth of the limit cycles.

The two bases therefore play complementary roles. The normal-mode basis is the physical language: it labels the lobe as antisymmetric and makes the sideband instability interpretable as a hybridization of normal-mode sidebands. The bare basis is the solvability language: it exposes the fact that resonator $2$ is not directly pumped, so its stationary equation can be inverted exactly. We will move between these descriptions several times. This is not a change of model, but a way of using the most transparent coordinates for each part of the argument.

Because $J<0$, equivalently $g<0$, the antisymmetric normal mode is the lower-frequency mode. In the notation of the normal-mode slow flow equations [cf.~Eqs.~\eqref{eq:nm_aa}--\eqref{eq:nm_as}],
\begin{equation}
  \Delta_a=\delta+g\,,
  \qquad
  \Delta_s=\delta-g\,,
\end{equation}
so the antisymmetric lobe is obtained by tuning $\Delta_a\simeq0$ and has $\Delta_s>0$. This is the correct physical language for the lobe. The most economical analytic calculation, however, uses the bare basis as a solvability gauge, because resonator $2$ is passive. We solve the stationary states in that basis and then translate the result back into the normal-mode language.

\subsection{Inversion symmetry}

Both the bare-basis slow flow~\eqref{eq:model_a1}--\eqref{eq:model_a2} and the normal-mode slow flow~\eqref{eq:nm_aa}--\eqref{eq:nm_as} are odd under simultaneous inversion of all slow amplitudes. In the normal-mode basis,
\begin{equation}
  (a_a,a_s)\longmapsto(-a_a,-a_s)\,.
  \label{eq:inv_sym}
\end{equation}
This is the discrete $\mathds{Z}_2$ symmetry inherited from the phase symmetry of the parametrically driven oscillator: shifting the carrier phase by $\pi$ flips the sign of every slow quadrature while leaving the equations invariant. It has two consequences that organize the rest of the analysis. Every nonzero stationary solution $(A_a,A_s)$ has a partner $(-A_a,-A_s)$, and every periodic orbit has an inversion-related partner. The two stationary phase states constructed below are precisely such a pair, as are the two limit cycles $I_+$ and $I_-$ ``born'' at the Hopf bifurcation. The equality of the two switching barriers has the same origin, provided that the noise source also respects the inversion symmetry.

\subsection{Exact stationary branch from the passive mode}

Throughout this section, capital letters denote stationary values of the corresponding slow envelopes: $A_1, A_2$ are the bare-mode stationary amplitudes ($a_1, a_2$ frozen on a phase state) and $A_s, A_a$ are the normal-mode stationary amplitudes. Bogoliubov fluctuations $\delta a_j$ around the stationary state are introduced in Sec.~\ref{sec:hopf}. We take the two bare detunings in Eqs.~\eqref{eq:model_a1}--\eqref{eq:model_a2} equal and write them as $\delta$. The stationary equations read
\begin{align}
  0&=(-\tilde\Gamma+i\delta)A_1+iK|A_1|^2A_1-i\mu A_1^*-igA_2\,,\label{eq:bare_stat_1}\\
  0&=(-\tilde\Gamma+i\delta)A_2+iK|A_2|^2A_2-igA_1\,.\label{eq:bare_stat_2}
\end{align}
The second equation contains no direct pump: the passive resonator is driven only through the coupling $g$ to resonator $1$. Equation~\eqref{eq:bare_stat_2} can therefore be inverted algebraically to express $A_2$ in terms of $A_1$, reducing the two coupled complex stationary equations to a single scalar problem whose roots are dubbed ``stationary branches''.

We define the passive nonlinear detuning
\begin{equation}
  \vartheta_2=\delta+K|A_2|^2 \,.
\end{equation}
This parameter is the detuning actually seen by the passive resonator after its Duffing frequency shift. Moving along the stationary branch is equivalent to moving the passive resonator through this nonlinear detuning. The pump does not directly set the shape of the stationary state. It selects which value of $\vartheta_2$ is realized.

Then Eq.~\eqref{eq:bare_stat_2} reveals the exact response relation
\begin{equation}
  \frac{A_2}{A_1}=\frac{g}{\vartheta_2+i\tilde\Gamma}\,.
  \label{eq:passive_ratio}
\end{equation}
The entire nonzero stationary branch is therefore parametrized by the single real variable $\vartheta_2$.

The amplitudes follow immediately:
\begin{equation}
  |A_2|^2=\frac{\vartheta_2-\delta}{K}\,,
  \qquad
  |A_1|^2=\frac{\vartheta_2^2+\tilde\Gamma^2}{g^2}\frac{\vartheta_2-\delta}{K}\,.
  \label{eq:bare_amplitudes_vartheta}
\end{equation}
The physical branch (where the phase-space coordinates are real) is selected by $({\vartheta_2-\delta})/K>0$. The corresponding nonlinear detuning of the driven resonator is
\begin{equation}
  \vartheta_1(\vartheta_2)=\delta+K|A_1|^2
  =\delta+\frac{\vartheta_2^2+\tilde\Gamma^2}{g^2}(\vartheta_2-\delta)\,.
  \label{eq:vartheta1}
\end{equation}
Here $\vartheta_1$ is the Kerr-shifted detuning of the driven resonator (resonator $1$) and $\vartheta_2$ that of the passive resonator (resonator $2$); both reduce to the bare detuning $\delta$ in the linear limit $|A_j|\to0$. The branch is most naturally parametrized by $\vartheta_2$ because the passive equation~\eqref{eq:bare_stat_2} contains no direct pump.
Substituting Eq.~\eqref{eq:passive_ratio} into the driven stationary equation results in 
\begin{equation}
  \left[-\tilde\Gamma_{\rm eff}(\vartheta_2)+i\Theta_{\rm eff}(\vartheta_2)\right]A_1=i\mu A_1^*\,,
\end{equation}
where
\begin{align}
  \tilde\Gamma_{\rm eff}(\vartheta_2)=\tilde\Gamma\left(1+\frac{g^2}{\vartheta_2^2+\tilde\Gamma^2}\right)\,,\quad
  \rm{and}\quad
  \Theta_{\rm eff}(\vartheta_2)=\vartheta_1(\vartheta_2)-\frac{g^2\vartheta_2}{\vartheta_2^2+\tilde\Gamma^2}\,.
  \label{eq:effective_susceptibility}
\end{align}
Thus, the passive resonator appears as a self-energy that dresses both the damping and the detuning of the driven mode. Taking the modulus yields the exact pump curve
\begin{equation}
  \mu^2(\vartheta_2)=
  \tilde\Gamma_{\rm eff}(\vartheta_2)^2+\Theta_{\rm eff}(\vartheta_2)^2 \,.
  \label{eq:mu_vartheta}
\end{equation}
For a fixed pump, the nonzero stationary states are the real roots of $\mu^2=\mu^2(\vartheta_2)$ satisfying $({\vartheta_2-\delta})/K>0$. The origin, $A_1=A_2=0$, is a separate stationary solution.

% This last sentence is the practical algorithm. Choose the pump. Solve the scalar equation $\mu^2=\mu^2(\vartheta_2)$. Each physical root gives $|A_2|$ and $|A_1|$ from Eq.~\eqref{eq:bare_amplitudes_vartheta}, the relative amplitude $A_2/A_1$ from Eq.~\eqref{eq:passive_ratio}, and the absolute phase from the equation below. Thus the finite-amplitude stationary problem has been reduced from two coupled complex equations to one real algebraic root selection.

The squared modulus in Eq.~\eqref{eq:mu_vartheta} discards phase information, but the Bogoliubov matrix below depends on $A_j^2$, not only on $|A_j|^2$. The missing phase is supplied by the driven stationary equation. If $A_1=|A_1|e^{i\phi_1}$, then
\begin{equation}
  e^{-2i\phi_1}=\frac{-\tilde\Gamma_{\rm eff}(\vartheta_2)+i\Theta_{\rm eff}(\vartheta_2)}{i\mu} \,.
  \label{eq:bare_phase}
\end{equation}
Together, Eqs.~\eqref{eq:passive_ratio}, \eqref{eq:bare_amplitudes_vartheta}, \eqref{eq:mu_vartheta}, and \eqref{eq:bare_phase} solve the damped stationary states analytically.

\subsection{Normal-mode geometry of the solved branch}

To physically interpret the stationary state along stationary branch, we quantify its normal-mode composition via the complex ratio $z = A_s/A_a$. The proportionality $A_1\propto A_s+A_a$ and $A_2\propto A_s-A_a$ between the bare and normal-mode amplitudes implies that the ratio reads
\begin{equation}
  z=\frac{A_s}{A_a}=\frac{A_1+A_2}{A_1-A_2}=\frac{\vartheta_2+i\tilde\Gamma+g}{\vartheta_2+i\tilde\Gamma-g}\,,
  \label{eq:z_vartheta}
\end{equation}
where we used Eq.~\eqref{eq:passive_ratio}.
Thus, the normal-mode branch is not a generic curve. It is a circle in the complex $z$ plane:
\begin{equation}
  (\operatorname{Re}z-1)^2+\left(\operatorname{Im}z+\frac{g}{\tilde\Gamma}\right)^2=\left(\frac{g}{\tilde\Gamma}\right)^2 \,.
  \label{eq:z_circle}
\end{equation}
For the antisymmetric lobe $g<0$, this circle has its center at $(1,-g/\tilde\Gamma)$ and radius $|g|/\tilde\Gamma$.

The circle is a useful geometric way of reading the branch. In the absence of damping, the response of the passive resonator would be in phase or antiphase with the driven one, and the corresponding part of the $z$-branch would be real. Damping gives the passive oscillator a phase lag. In the normal-mode ratio~\eqref{eq:z_vartheta}, this lag appears as an imaginary part of $z$. The imaginary part is therefore not an additional branch parameter. It is fixed by dissipation once $\vartheta_2$ has been chosen.

The same ratio also tells us how the phase state is distributed over the two physical resonators:
\begin{equation}
  \frac{A_2}{A_1}=\frac{z-1}{z+1}\,.
  \label{eq:bare_ratio_from_z}
\end{equation}
Near the antisymmetric lobe threshold, $z\simeq0$ and $A_2/A_1\simeq-1$. The state is an antisymmetric normal mode. As $|\vartheta_2|$ grows, $A_2/A_1\to0$, and the state becomes localized on the directly driven resonator. The stationary branch is therefore not merely an antisymmetric state with a small symmetric correction. It is a continuous migration from a collective antisymmetric mode to a driven-resonator-localized state.

\subsection{Threshold and critical eigenvector}

For the stationary branch at the origin, the Kerr shifts vanish. In the general bare-basis problem, the passive resonator dresses the driven resonator susceptibility, and the parametric threshold is
\begin{equation}
  \mu_{\rm th}^2=\left|\tilde\Gamma-i\delta_1+\frac{g^2}{\tilde\Gamma-i\delta_2}\right|^2\,.
  \label{eq:threshold_general}
\end{equation}
In the decoupled limit $g\to0$, this reduces to the familiar single-oscillator result $\mu_{\rm th}^2=\tilde\Gamma^2+\delta_1^2$. For finite coupling, the self-energy $g^2/(\tilde\Gamma-i\delta_2)$ shifts the effective damping and detuning of the driven resonator. This is the linear version of the finite-amplitude self-energy in Eq.~\eqref{eq:effective_susceptibility}.

For the antisymmetric lobe, the same result follows directly from the stationary branch by setting $\vartheta_2=\delta$. Since $\Delta_a=\delta+g=0$, we have $\delta=-g$, and Eq.~\eqref{eq:mu_vartheta} gives
\begin{equation}
  \mu_{{\rm th},a}^2=\frac{\tilde\Gamma^2(\tilde\Gamma^2+4g^2)}{\tilde\Gamma^2+g^2}\,.
  \label{eq:threshold_a}
\end{equation}
The threshold is therefore not an extra calculation separate from the nonlinear branch. It is the endpoint of the same construction. At threshold, the nonlinear shifts vanish, $\vartheta_2=\delta$, and the exact finite-amplitude pump curve reduces to the linear susceptibility condition. The corresponding symmetric-lobe threshold $\mu_{{\rm th},s}$ follows by setting $\Delta_s=0$ instead, which gives $\delta=g$ and the same expression with $g\to-g$; since the formula is even in $g$, the two lobe thresholds coincide.

The composition of the emerging phase state at threshold is equally direct. The destabilizing eigenmode of the linearized dynamics at $\mu=\mu_{{\rm th},a}$ has a fixed normal-mode ratio $z_{{\rm th},a}=A_s/A_a$, obtained by evaluating Eq.~\eqref{eq:z_vartheta} at $\vartheta_2=\delta=-g$:
\begin{equation}
  z_{{\rm th},a}=\frac{\tilde\Gamma}{\tilde\Gamma+2ig}\,.
  \label{eq:z_threshold_a}
\end{equation}
Thus, the critical state is mostly antisymmetric only when $|g|\gg\tilde\Gamma$. For finite damping, the symmetric normal-mode component is already present at threshold, with magnitude $|z_{{\rm th},a}|=\tilde\Gamma/\sqrt{\tilde\Gamma^2+4g^2}$.

\section{Linear response, Hopf bifurcation, and the birth of limit cycles}
\label{sec:hopf}

\subsection{Bogoliubov matrix on the solved branch}

The stationary branch has now been solved. The next question is when it loses stability. The bare basis is the most compact for this purpose because the nonlinearities are local. Writing $a_j = A_j(\vartheta_2) + \delta a_j$ and linearizing Eqs.~\eqref{eq:model_a1}--\eqref{eq:model_a2}, we collect the fluctuations into the Nambu vector $Y=(\delta a_1,\delta a_1^*,\delta a_2,\delta a_2^*)^T$. The linearized dynamics is $\dot Y = LY$ with
\begin{equation}
  L=\begin{pmatrix}
    id_1-\tilde\Gamma & U_1 & -ig & 0\\
    U_1^* & -id_1-\tilde\Gamma & 0 & ig\\
    -ig & 0 & id_2-\tilde\Gamma & U_2\\
    0 & ig & U_2^* & -id_2-\tilde\Gamma
  \end{pmatrix}\,.
  \label{eq:L_bare}
\end{equation}
The diagonal carries uniform damping $-\tilde\Gamma$ together with the amplitude-shifted detunings $d_j = 2\vartheta_j - \delta$; the off-diagonal $\pm ig$ is the ordinary hopping that hybridizes the two resonators. The anomalous (squeezing) entries $U_1 = iKA_1^2 - i\mu$ and $U_2 = iKA_2^2$ come from the curvature of the Kerr nonlinearity around the stationary state, with the external pump $\mu$ entering only through $U_1$. All entries are algebraic functions of the single branch parameter $\vartheta_2$; in particular $U_1$ and $U_2$ depend on $A_j^2$ rather than only $|A_j|^2$, which is why the phase reconstruction in Eq.~\eqref{eq:bare_phase} was needed. The Hopf bifurcation is a sideband instability of this Bogoliubov problem as the stationary state moves along the branch.

Because the damping is uniform, the spectrum of $L$ is symmetric about $-\tilde\Gamma$ and its characteristic polynomial is biquadratic in the shifted eigenvalue $q + \tilde\Gamma$:
\begin{equation}
  \det(q\mathds{1}-L)=(q+\tilde\Gamma)^4+\sigma(\vartheta_2)(q+\tilde\Gamma)^2+\pi(\vartheta_2)\,,
  \label{eq:shifted_biquadratic}
\end{equation}
with invariants $\sigma(\vartheta_2) = -\tfrac{1}{2}\operatorname{tr}(L+\tilde\Gamma\mathds{1})^2$ and $\pi(\vartheta_2) = \det(L+\tilde\Gamma\mathds{1})$. Instead of tracking a generic quartic along a numerically continued branch, the bifurcations become scalar algebraic conditions on the two invariants $\sigma$ and $\pi$ evaluated on the exact stationary branch. The invariants also have a direct spectral meaning: in the undamped (Hamiltonian) limit, the two sideband-frequency squares have sum $\sigma$ and product $\pi$, so the curve $\vartheta_2\mapsto(\sigma(\vartheta_2),\pi(\vartheta_2))$ tracks how the sideband pairs move as the stationary phase state hybridizes, and static instabilities, sideband collisions, and Hopf bifurcations appear as distinct boundaries in this plane.

% For readers who want the Routh--Hurwitz form explicitly, Eq.~\eqref{eq:shifted_biquadratic} is equivalent to
% \begin{equation}
%   q^4+c_1q^3+c_2q^2+c_3q+c_4,
% \end{equation}
% with
% \begin{equation}
%   c_1=4\tilde\Gamma,
%   \qquad
%   c_2=6\tilde\Gamma^2+a,
%   \qquad
%   c_3=2\tilde\Gamma(a+2\tilde\Gamma^2),
%   \qquad
%   c_4=\tilde\Gamma^4+a\tilde\Gamma^2+b.
%   \label{eq:quartic_coefficients}
% \end{equation}
% Thus one may also define $R(\vartheta_2)=\sigma(\vartheta_2)+2\tilde\Gamma^2$, for which $c_3=2\tilde\Gamma R$. This connects the invariant notation to the more standard quartic stability test.

\subsection{Static and Hopf bifurcations}

The biquadratic structure of Eq.~\eqref{eq:shifted_biquadratic} reduces both bifurcations of the stationary phase state to scalar conditions on the invariants $\sigma$ and $\pi$. A static instability ($q=0$) implies
\begin{equation}
  \pi(\vartheta_2)=-\tilde\Gamma^2\left[\sigma(\vartheta_2)+\tilde\Gamma^2\right]\,,
  \label{eq:static_bifurcation}
\end{equation}
which is the linear-response form of a fold condition along the one-parameter branch, away from symmetry-enforced points such as the origin threshold. A Hopf bifurcation ($q=\pm i\Omega_H$) requires
\begin{equation}
  \pi(\vartheta_H)=\tfrac{1}{4}\left[\sigma(\vartheta_H)+4\tilde\Gamma^2\right]^2\,,
  \qquad
  \Omega_H^2=\tilde\Gamma^2+\tfrac{1}{2}\sigma(\vartheta_H)\,,
  \label{eq:hopf}
\end{equation}
with corresponding pump value $\mu_H=\mu(\vartheta_H)$ from Eq.~\eqref{eq:mu_vartheta}. Equivalently, the Hopf point is the physical root of $\Psi(\vartheta_2)\equiv\pi(\vartheta_2)-\tfrac{1}{4}[\sigma(\vartheta_2)+4\tilde\Gamma^2]^2$ on a branch stable up to that boundary.

The full prediction is purely algebraic: pick $\vartheta_2$, obtain $A_1, A_2$ from Eq.~\eqref{eq:bare_amplitudes_vartheta}, evaluate $\sigma$ and $\pi$ on $L(\vartheta_2)$, then solve Eq.~\eqref{eq:hopf} for $\vartheta_H$ to read off $\mu_H$ and $\Omega_H$. The Hopf is reached when a pair of opposite-signature sidebands couples strongly enough that one damping rate vanishes.

\subsection{Sideband attraction and exceptional points}

The eigenvalues of $L$ are sidebands of the stationary phase state. In the undamped limit, a sideband collision occurs when the two roots of the biquadratic in Eq.~\eqref{eq:shifted_biquadratic} coalesce, $\sigma(\vartheta_2)^2 = 4\pi(\vartheta_2)$. This is the exceptional-point or Krein-collision condition of the undamped Bogoliubov problem; finite damping shifts the actual Hopf beyond it, following the sequence frequency attraction $\to$ exceptional point $\to$ linewidth splitting $\to$ Hopf, reached only once the linewidth split from the opposite-signature collision overcomes the uniform damping $\tilde\Gamma$.

Whether two colliding sidebands repel or attract is set by the Nambu (Krein) signature $\kappa_v = v^\dagger\Sigma_z v$ with $\Sigma_z = \operatorname{diag}(1,-1,1,-1)$, evaluated on the normalized eigenvector $v$. Same-signature sidebands repel and form avoided crossings; opposite-signature sidebands attract, lock in frequency, and split in damping~\cite{Krein1950,mackay2020stability,kapitula2013spectral,Soriente_2020,SorienteDistinctive2021,dumont2024hamiltonian}. Here, the attracting sidebands belong to the hybridizing phase state described by Eq.~\eqref{eq:z_vartheta}, so the Hopf is the sideband instability of a phase state migrating from the antisymmetric normal mode toward the driven resonator,.

The self-energy structure is transparent in the simplified limit $U_2\approx 0$. Eliminating the passive sidebands gives
\begin{equation}
  \left[q+\tilde\Gamma-id_1+\frac{g^2}{q+\tilde\Gamma-id_2}\right]\left[q+\tilde\Gamma+id_1+\frac{g^2}{q+\tilde\Gamma+id_2}\right]-|U_1|^2=0\,,
  \label{eq:factorized_dispersion}
\end{equation}
whose bracketed factors are the positive- and negative-frequency susceptibilities of the driven phase state, each dressed by a passive-resonator self-energy, mixed by the anomalous term $|U_1|^2$. The passive resonator therefore supplies an additional sideband that can collide with a sideband of the driven state, the relative Krein signature deciding between avoided crossing and level attraction.

Near a two-sideband collision the local form is
\begin{equation}
  q_\pm=-\tilde\Gamma+i\bar\Omega\pm\sqrt{|\kappa|^2-\left(\tfrac{\Delta}{2}\right)^2}\,,
  \label{eq:level_attraction}
\end{equation}
where $\Delta$ is the detuning between the uncoupled sidebands, $\kappa$ their effective coupling, and $\bar\Omega$ their mean frequency. The exceptional point sits at $|\kappa|=|\Delta|/2$ and the Hopf at $|\kappa|^2 = \tilde\Gamma^2 + (\Delta/2)^2$, the local two-mode version of the invariant condition~\eqref{eq:hopf}. This is the same level-attraction mechanism that organizes sideband instabilities via internal resonance in nonlinear mechanical oscillators~\cite{fu2025sideband}; the coupled-Kerr parametric landscape is discussed more broadly in Ref.~\cite{Ameye2025Parametric}.

\subsection{From phase states to limit cycles}
\label{sec:phase_to_cycles}

We can now assemble the deterministic picture. Since $J<0$, the experimentally relevant lobe is the antisymmetric normal-mode lobe. The origin first loses stability at the threshold in Eq.~\eqref{eq:threshold_a}, and the emerging nonzero phase state has the critical ratio in Eq.~\eqref{eq:z_threshold_a}. The full finite-amplitude phase state is then generated by the passive-mode parameter $\vartheta_2$ and mapped to the normal-mode ratio $z(\vartheta_2)$ in Eq.~\eqref{eq:z_vartheta}. As the pump is increased, the state continuously hybridizes and becomes increasingly localized on the driven resonator.

The stability of this solved branch is governed by the two invariants $\sigma(\vartheta_2)$ and $\pi(\vartheta_2)$. At the Hopf boundary~\eqref{eq:hopf}, an opposite-signature sideband collision produces a damping split large enough to overcome $\tilde\Gamma$. Provided the Hopf bifurcation is supercritical, each stationary phase state is replaced by a stable limit cycle. The two resulting cycles,
\begin{equation}
  I_+ \quad \text{and} \quad I_-\,,
\end{equation}
are related by the inversion symmetry~\eqref{eq:inv_sym}. The linear Hopf calculation fixes their onset and sideband frequency. The nonlinear Hopf normal form controls their initial amplitude and frequency pulling, but the global switching barrier between them is a separate large-deviation problem.

Our deterministic analysis has now accomplished everything that a local bifurcation calculation can do. It explains why the two attractors are born in symmetry-related pairs, how the stationary phase states are organized, and how the Hopf onset is predicted analytically. It does not, however, determine how a noisy trajectory moves from one attracting cycle to the other. That question depends on the global basin geometry, not only on the local Hopf neighborhood.

With the limit cycles $I_\pm$ in hand, the deterministic skeleton is complete. The switching barrier is the height of the action landscape between basins, and the optimal transition path is a global object that must find the cheapest route to the separatrix. These are precisely the quantities addressed by Freidlin-Wentzell large-deviation theory.

\bibliography{aipsamp}% Produces the bibliography via BibTeX.